\documentclass[12pt]{article}

\usepackage{latexsym, amssymb, amsmath, amscd, amsfonts, epsfig, graphicx, colordvi,verbatim,ifpdf}
\usepackage{amsfonts, amsmath, amssymb}
\usepackage{amssymb,amsfonts,amsmath,latexsym,epsfig,cite, psfrag,eepic,color}
\usepackage{amscd,graphics}
\usepackage{latexsym, amssymb,  amsmath,amscd, amsfonts, epsfig, graphicx, colordvi,amsthm}

\usepackage{graphicx}
\usepackage{color}
\usepackage{ifpdf}
\usepackage{fancybox}

\newtheorem{thm}{Theorem}[section]

\newtheorem{core}[thm]{Corollary}

\def\pf{\noindent{\it Proof.} }
\def\qed{\nopagebreak\hfill{\rule{4pt}{7pt}}
\medbreak}

\numberwithin{equation}{section}

\def\qed{\nopagebreak\hfill{\rule{4pt}{7pt}}
\medbreak}

\def \c {NC}

\def \f {F}

\newlength{\boxedparwidth}
  {\begin{center} \begin{tabular}{|@{\hspace{.315in}}c@{\hspace{.15in}}|}
                  \hline \\ \begin{minipage}[t]{\boxedparwidth}
                  \setlength{\parindent}{.25in}}%
  {\end{minipage} \\ \\ \hline \end{tabular} \end{center}}

\begin{document}
\begin{center}

{\Large \bf $k$-Marked Dyson Symbols and

Congruences for Moments of Cranks }
\end{center}

\vskip 5mm

\begin{center}
{  William Y.C. Chen}$^{1}$,    {  Kathy Q. Ji}$^{2}$
and {   Erin  Y.Y. Shen}$^{3}$ \vskip 2mm

    $^{1,2,3}$Center for Combinatorics, LPMC-TJKLC\\
   Nankai University, Tianjin 300071, P. R. China\\[6pt]
   $^{1}$Center for Applied Mathematics\\
Tianjin University,  Tianjin 300072, P. R. China\\[6pt]

   \vskip 2mm

    $^1$chen@nankai.edu.cn, $^2$ji@nankai.edu.cn,  $^3$shenyiying@mail.nankai.edu.cn
\end{center}

\vskip 6mm \noindent {\bf Abstract.}  By introducing
$k$-marked Durfee symbols, Andrews found a
 combinatorial interpretation of   $2k$-th symmetrized  moment $\eta_{2k}(n)$ of ranks of partitions of $n$. Recently, Garvan introduced the $2k$-th symmetrized moment $\mu_{2k}(n)$ of cranks of partitions of $n$ in the study of the higher-order spt-function $spt_k(n)$.    In this paper, we give a combinatorial interpretation of $\mu_{2k}(n)$. We introduce  $k$-marked Dyson symbols  based on  a  representation of ordinary partitions given by Dyson, and we show  that   $\mu_{2k}(n)$ equals the number of $(k+1)$-marked Dyson symbols of $n$. We then  introduce  the full crank of a $k$-marked Dyson symbol and show that there exist an infinite family of congruences for the full crank function of $k$-marked Dyson symbols  which implies that  for fixed prime $p\geq 5$ and positive integers $r$ and $k\leq (p-1)/2$, there exist infinitely many non-nested arithmetic progressions $An+B$ such that $\mu_{2k}(An+B)\equiv 0\pmod{p^r}.$

\section{Introduction}

Dyson's rank \cite{Dyson-1944}  and   the Andrews-Garvan-Dyson  crank \cite{Andrews-Garvan-1988} are two fundamental statistics in the theory of partitions.  For a partition $\lambda=(\lambda_1,\lambda_2,\ldots,\lambda_\ell)$, the rank of
$\lambda$, denoted $r(\lambda)$, is   the largest part of $\lambda$ minus the number of parts.  The crank $c(\lambda)$     is defined  by
 \[c(\lambda)=\left\{
\begin{array}{ll}
\lambda_1,\ \ &\text{ if}\ n_1(\lambda)=0,\\[10pt]
 \mu(\lambda)-n_1(\lambda), \ &\text{ if }\ n_1(\lambda)>0,
\end{array}\right.
\]
where   $n_1(\lambda)$ is  the number of ones in $\lambda$ and $\mu(\lambda)$ is  the number of parts  larger than $n_1(\lambda)$.

Andrews \cite{Andrews-07-a} introduced  the symmetrized moments $\eta_{2k}(n)$ of ranks of partitions of $n$  given by
\begin{align}\label{symRankMom}
\eta_{k}(n)=\sum_{m=-\infty}^{+\infty}{m+\lfloor\frac{k-1}{2}\rfloor \choose k}N(m,n),
\end{align}
where $N(m,n)$ is the number of partitions of
$n$ with rank $m$.

In view of the symmetry $N(-m,n)=N(m,n)$,  we have $\eta_{2k+1}(n)=0$.  As for the even symmetrized
moments $\eta_{2k}(n)$, Andrews \cite{Andrews-07-a}  showed that for fixed $k\geq 1$,
  $\eta_{2k}(n)$  is equal to the number of $(k+1)$-marked Durfee symbols of $n$.   Kursungoz \cite{Kursungoz-2011} and Ji \cite{Ji-2011} provided the alternative proof of this result respectively. Bringmann, Lovejoy and Osburn \cite{Bringmann-2010} defined two-parameter generalization of $\eta_{2k}(n)$ and $k$-marked Durfee symbols. In \cite{Andrews-07-a},  Andrews also introduced the full rank of a $k$-marked Durfee symbol and defined  the full rank function $NF_k(r,t;n)$ to be  the number of $k$-marked Durfee symbols of $n$ with full rank congruent to $r$ modulo $t$.

The full rank function $NF_k(r,t;n)$  have been extensively studied and they posses many congruence properties, see for example, \cite{Bringmann-2008, Bringmann-Garvan-Mahlburg-2009, Bringmann-2010, Bringmann-Kane-2012, Keith-2009}.    Recently,   Bringmann, Garvan and Mahlburg \cite{Bringmann-Garvan-Mahlburg-2009} used the automorphic properties of the generating functions of $NF_k(r,t;n)$ to prove the existence of infinitely many congruences for $NF_k(r,t;n)$. More precisely, for given positive integers $j$, $k\geq 3$, odd positive integer $t$, and prime $Q$ not divisible by $6t$, there exist infinitely many arithmetic progressions $An+B$ such that for every $0\leq r<t$, we have
\begin{equation}\label{nfk}
NF_k(r,t;An+B)\equiv 0\pmod{Q^j}.
\end{equation}
Since
\[
\eta_{2k}(n)=\sum_{r=0}^{t-1}NF_{k+1}(r,t;n),
\] by (\ref{nfk}),  we see
 that there exist  an  infinite family of congruences for $\eta_{2k}(n)$, namely, for given positive integers $k$ and $j$, prime $Q>3$, there exist infinitely many non-nested arithmetic progressions $An+B$ such that
\[\eta_{2k}(An+B)\equiv 0 \pmod{Q^j}.\]

Analogous to the symmetrized moments $\eta_{k}(n)$ of ranks,  Garvan \cite{Garvan-2011} introduced the $k$-th symmetrized moments $\mu_k(n)$ of cranks of partitions of $n$  in the study of the higher-order spt-function $spt_k(n)$.  To be more specific,
\begin{align}\label{symRankMom}
\mu_{k}(n)=\sum_{m=-\infty}^{+\infty}{m+\lfloor\frac{k-1}{2}\rfloor \choose k}M(m,n),
\end{align}
where $M(m,n)$ denotes the number of partitions of $n$ with crank $m$ for $n> 1$. For $n=1$ and $m\neq -1,0,1$, we set $M(m,1)=0$; otherwise,  we define
\[M(-1,1)=1,\  M(0,1)=-1, \ M(1,1)=1.\]
 It is clear that $\mu_{2k+1}(n)=0$, since $M(m,n)=M(-m,n)$.

In this paper,  we give a combinatorial interpretation of $\mu_{2k}(n)$.
 We first introduce the notion of $k$-marked Dyson symbols based on  a  representation for ordinary partitions given by Dyson \cite{Dyson-1944}. We  show that for fixed $k\geq 1$, $\mu_{2k}(n)$ equals the number of $(k+1)$-marked Dyson symbols of $n$. Moreover, we define the full crank of a $k$-marked Dyson symbol and
define full crank function $NC_k(r,t;n)$ to be the number of $k$-marked Dyson symbols of $n$ with full crank congruent to $r$ modulo  $t$. We prove that for fixed prime $p\geq 5$ and positive integers $r$ and $k\leq (p+1)/2$,     there exists infinitely many non-nested arithmetic progressions $An+B$ such that for every $0\leq i\leq p^r-1$,
\begin{equation}\label{nck}
\c_k(i, p^r; An+B) \equiv 0 \pmod{p^r}.
\end{equation}
Note that
\[
\mu_{2k}(n)=\sum_{i=0}^{p^r-1}\c_{k+1}(i,p^r;n),
\]
so that from \eqref{nck} we can deduce  that   there exist an  infinite family of congruences for $\mu_{2k}(n)$, that is, for  fixed prime $p\geq 5$,  positive integers $r$ and $k\leq (p-1)/2$, there exist infinitely many non-nested arithmetic progressions $An+B$ such that
\begin{equation*}
\mu_{2k}(An+B)\equiv 0 \pmod{p^r}.
\end{equation*}

\section{Dyson symbols and $k$-marked Dyson symbols}

In this section, we introduce the notion of $k$-marked Dyson symbols.
A $1$-marked Dyson symbol is called a Dyson symbol,
which is a representation of a   partition introduced by  Dyson \cite{Dyson-1989}.  For $1\leq i\leq k$, we  define the $i$-th crank of a $k$-marked Dyson symbol. Moreover, we define the function $F_k(m_1,m_2,\ldots,m_{k};n)$ to be the number of $k$-marked Dyson symbol of $n$ with the $i$-th  crank equal to  $m_i$ for $1\leq i \leq k$.
 The following theorem shows that the number of $k$-marked Dyson symbols of $n$
 can be expressed in terms of the number of Dyson symbols of $n$.

 \begin{thm} \label{main} For fixed integers $m_1,m_2,\ldots, m_k$, we have
 \begin{equation}\label{main-e}
\f_k(m_1,\ldots,m_k;n)=\sum_{t_1,\ldots,\,t_{k-1}=0}^{+\infty}
\f_1\left(\sum_{i=1}^k|m_i|+2\sum_{i=1}^{k-1}t_i+k-1;n\right).
\end{equation}
\end{thm}

For a   partition $\lambda=(\lambda_1,\ldots,\lambda_\ell)$, let   $\ell(\lambda)$ denote the number of parts of $\lambda$ and $|\lambda|$ denote the sum of parts of $\lambda$.  A Dyson symbol of $n$ is a pair of restricted partitions $(\alpha,\beta)$ satisfying the following conditions:
\begin{itemize}
\item[(1)]If  $\ell(\alpha)=0$, then $\beta_1=\beta_2$;
\item[(2)]If $\ell(\alpha)=1$, then $\alpha_1=1$;
\item[(3)]If $\ell(\alpha)>1$, then $\alpha_1=\alpha_2$;
\item[(4)]$n=|\alpha|+|\beta|+\ell(\alpha) \ell(\beta)$.
\end{itemize}
When we display a Dyson symbol, we shall put $\alpha$ on the top of $\beta$
in the form of a Durfee symbol \cite{Andrews-07-a} or a Frobenius partition \cite{Andrews-1984}.

For example, there are $5$ Dyson symbols   of $4$:
\[\left(\begin{array}{ccc}
&\\
2&2
\end{array}
\right),\ \left(\begin{array}{cccc}
&&&\\
1&1&1&1
\end{array}
\right),\ \left(\begin{array}{c}
1\\
2
\end{array}
\right),\ \ \left(\begin{array}{ccc}
2&2 \\
&
\end{array}
\right),\ \left(\begin{array}{cccc}
1&1&1&1 \\
&&&
\end{array}
\right).\]

\begin{thm}[Dyson] \label{dy-op}There is a bijection $\Omega$ between the set of partitions of $n$ and the set of Dyson symbols of $n$.
\end{thm}

For completeness, we give a proof of the above theorem.

\noindent {\it Proof of Theorem \ref{dy-op}$\colon$ }
 Let  $\lambda=(\lambda_1,\ldots, \lambda_\ell)$ be a partition of $n$.
  A Dyson symbol $(\alpha,\beta)$ of $n$ can be constructed via the
  following procedure.   There are two cases.

\noindent Case 1: One is not a part of $\lambda$.
 We set $\alpha=\emptyset$ and $\beta=\lambda'$.

\noindent Case 2:   One is a part of $\lambda$. Assume that
  one occurs $M$ times in $\lambda$. We decompose the Ferrers diagram of
   $\lambda$ into three
  blocks as illustrated in  Figure 2.1, where $N$ is the number of parts of $\lambda$ that are greater than $M$. In this case, we see that   $\lambda=(\lambda_1,\ldots,\lambda_N,\lambda_{N+1}, \ldots, \lambda_s, 1^M)$, where $\lambda_N>M$,  $\lambda_{N+1}\leq M$ and $1^M$ means $M$ occurrences of $1$. Then
  remove all parts equal to one  from $\lambda$ and insert a new part $M$,
  so that  we get a   partition  $\mu=(\lambda_1,\ldots,\lambda_N, M, \lambda_{N+1}, \ldots, \lambda_s)$ as shown in Figure 2.2.

\input{decomp-n.TpX}

Now the partitions $\alpha$ and $\beta$ can be obtained from  $\mu$. First, let $\beta=(\lambda_1-M, \lambda_2-M, \ldots, \lambda_N-M)$,
and let $\nu=(M, \lambda_{N+1}, \ldots, \lambda_s)$.
Then we get $\alpha=(\nu'_1, \nu'_2,\ldots, \nu'_M)$,
where $\nu'$  the conjugate of $\nu$, see Figure 2.2.

\input{decomp-2-n.TpX}

It is easy to verify that $(\alpha,\beta)$ is a Dyson symbol of $n$ and the above procedure is
 reversible, and hence the proof is complete.  \qed

For a Dyson symbol $(\alpha,\,\beta)$, Dyson \cite{Dyson-1989} considered  the difference between the number of parts of $\alpha$ and
$\beta$, which we call  the crank of $(\alpha,\,\beta)$.  Let $\f_1(m;n)$  denote the
number of  Dyson symbols of $n$ with crank $m$. Dyson \cite{Dyson-1989} observed the following relation based on the construction   in Theorem \ref{dy-op}.

\begin{core}[Dyson]\label{core-dyson} For $n\geq 2$ and integer $m$,
\begin{equation}\label{mf}
M(-m,n)=F_1(m;n).
\end{equation}
\end{core}

A $k$-marked  Dyson symbol  is defined as the following array
 \[\eta=\left(\begin{array}{cccccccc}
\alpha^{(k)},&&\alpha^{(k-1)},&&\ldots,&&\alpha^{(1)}\\[2pt]
&p_{k-1},&&p_{k-2},&\cdots&p_1,\\[2pt]
\beta^{(k)},&&\beta^{(k-1)},&&\ldots,&&\beta^{(1)}
\end{array}
\right),\]
consisting of  $k$ pairs of
partitions $(\alpha^{(i)},\beta^{(i)})$  and a partition $p=(p_{k-1},p_{k-2},\ldots, p_{0})$
 subject to the following conditions:
\begin{itemize}

\item[{\rm (1)}] The smallest part of $p$ equals $1$, that is, $p_{k-1} \geq \cdots \geq p_1\geq p_{0}=1$.

\item[{\rm (2)}]For $1\leq i\leq k-1$,  each part of $\alpha^{(i)}$ and $\beta^{(i)}$ is between $p_{i-1}$ and $p_i$, namely,
    \[p_i\geq \alpha^{(i)}_{1}\geq \alpha^{(i)}_{2}\geq \cdots \geq \alpha^{(i)}_{\ell}\geq p_{i-1} \quad \text{and} \quad  p_i\geq \beta^{(i)}_{1}\geq \beta^{(i)}_{2}\geq \cdots \geq \beta^{(i)}_{\ell}\geq p_{i-1}.\]

\item[{\rm (3)}] Each part of $\alpha^{(k)}$ and $\beta^{(k)}$ is
no  less  than $p_{k-1}$, namely,
 \[\alpha^{(k)}_1\geq \alpha^{(k)}_2 \geq \cdots \geq \alpha^{(k)}_\ell\geq p_{k-1} \quad \text{and} \quad \beta^{(k)}_1\geq \beta^{(k)}_2 \geq \cdots \geq \beta^{(k)}_\ell\geq p_{k-1}.\]

\item[{\rm (4)}]
If  $\ell(\alpha^{(k)})=1,$ then $\alpha^{(k)}_1=p_{k-1}${\rm ;}\\[3pt]
If $\ell(\alpha^{(k)})>1,$ then $\alpha^{(k)}_1=\alpha^{(k)}_2${\rm ;}\\[3pt]
If $\ell(\alpha^{(k)})=0$  and $\ell(\beta^{(k)})=1$, then
$\beta_1^{(k)}=p_{k-1}$;\\[3pt]
If $\ell(\alpha^{(k)})=0$ and $\ell(\beta^{(k)})\geq 2$, then  $\beta^{(k)}_1=\beta^{(k)}_2${\rm ;} \\[3pt]
If $\ell(\alpha^{(k)})=0$ and $\ell(\beta^{(k)})=0$, then $p_{k-1}=\max\{\alpha^{(k-1)}_1, \beta_1^{(k-1)}\}${\rm.}

\end{itemize}

For example, the array below
\begin{equation}\label{exam-eta}
\eta=\left(\begin{array}{cccccccc}
(5,5,4)&&(3,3,2)&&(1,1)\\[2pt]
&4&&2\\[2pt]
(4)&&(3,2,2)&&(2,1,1)
\end{array}\right)
\end{equation}
is a $3$-marked Dyson symbol.

We next define the weight of a $k$-marked Dyson symbol.
Recall that for a pair of partitions $(\alpha,\beta)$ with $\ell(\alpha)\geq \ell(\beta)$,
a balanced part $\beta_i$ of $\beta$ is defined recursively as follow.
If  the number of parts greater than
$\beta_{i}$ in $\alpha$   {\it is equal
to} the number of unbalanced parts   before $\beta_i$ in
$\beta$, that is, the number of unbalanced parts $\beta_j$ with $1\leq j<i$;
 otherwise, we call $\beta_i$ is an unbalanced part, see \cite[p.992]{Ji-2011}. We use
 $b(\alpha,\beta)$ to denote the number of balanced parts of $(\alpha,\beta)$.

 For example, for the pair of partitions
 \[\left(\begin{array}{c}\alpha\\ \beta
 \end{array}
 \right)=\left(\begin{array}{ccccc}
3&3&1&1\\
3&2&2&
\end{array}\right),\]   the first part  $3$ of $\beta$ is balanced, and the second part $2$ and the third part $2$ are unbalanced.  Therefore, $b(\alpha,\beta)=1.$

 We now define the $i$-th crank  and the $i$-th balanced number  of a $k$-marked Dsyon symbol.  Let
 \[\eta=\left(\begin{array}{cccccccc}
\alpha^{(k)},&&\alpha^{(k-1)},&&\ldots,&&\alpha^{(1)}\\[2pt]
&p_{k-1},&&p_{k-2},&\cdots&p_1\\[2pt]
\beta^{(k)},&&\beta^{(k-1)},&&\ldots,&&\beta^{(1)}
\end{array}
\right)\]
be a $k$-marked Dyson symbol. The pair of partitions
$(\alpha^{(i)},\beta^{(i)})$ is called the $i$-th vector of $\eta$.
For $1\leq i\leq k$, we define $c_i(\eta)$,  the $i$-th {\it crank}  of $\eta$,  to be the difference between the number of parts of $\alpha^{(i)}$ and $\beta^{(i)}$, that is, $c_i(\eta)=\ell(\alpha^{(i)})-\ell(\beta^{(i)})$.

 For $1\leq i<k$, we define  $b_i(\eta)$, the $i$-th balanced number of $\eta$ by
   \[b_i(\eta)= \begin{cases}b(\alpha^{(i)},\beta^{(i)}), & \text{if }\ell(\alpha^{(i)})\geq \ell(\beta^{(i)}),\\[5pt]
 b(\beta^{(i)}, \alpha^{(i)}), & \text{if }\ell(\alpha^{(i)})< \ell(\beta^{(i)}).
    \end{cases}\]
For $i=k$, we set $b_k(\eta)=0$.

For the $3$-marked Dyson symbol $\eta$ in \eqref{exam-eta},  we have $c_1(\eta)=-1,\, c_2(\eta)=0,\,c_3(\eta)=2$ and $b_1(\eta)=1,\,b_2(\eta)=1,\,b_3(\eta)=0$.

For $1\leq i\leq k$, we define $l_i(\eta)$, the $i$-th   large length   of $\eta$ by
\[l_i(\eta)= \begin{cases}\ell(\alpha^{(i)}), & \text{if }\ell(\alpha^{(i)})\geq \ell(\beta^{(i)}),\\[5pt]
 \ell(\beta^{(i)}), & \text{if }\ell(\alpha^{(i)})< \ell(\beta^{(i)}).
    \end{cases}\]
Similarly, we define the $i$-th    small  length  $s_i(\eta)$ of $\eta$ by
\[s_i(\eta)= \begin{cases}\ell(\beta^{(i)}), & \text{if }\ell(\alpha^{(i)})\geq \ell(\beta^{(i)}),\\[5pt]
 \ell(\alpha^{(i)}), & \text{if }\ell(\alpha^{(i)})< \ell(\beta^{(i)}).
    \end{cases}\]

The weight of $k$-marked Dyson symbol is defined by
\begin{equation}
|\eta|=\sum_{i=1}^k(|\alpha^{(i)}|+|\beta^{(i)}|)+\sum_{i=1}^{k-1}p_i
+(l(\eta)+D+k-1)(s(\eta)-D),
\end{equation}
where
\begin{equation}\label{eq-t}
 l(\eta)=\sum_{i=1}^k l_i(\eta), \quad   s(\eta)= \sum_{i=1}^k s_i(\eta), \quad \text{and} \quad D=\sum_{i=1}^kb_i(\eta).
 \end{equation}
For example, the weight of the $3$-marked Dyson symbol $\eta$ in  \eqref{exam-eta} equals $97$. 

For a $k$-marked Dyson symbol $\eta$, if the weight of $\eta$ equals $n$, we call   $\eta$  a $k$-marked Dyson symbol of $n$.
We can now define the function  $\f_k(m_1,\ldots,m_k;n)$ as the number of $k$-marked Dyson symbols of $n$ with the $i$-th crank equal to $m_i$ for $1\leq i \leq k$. Note that a $1$-marked Dyson symbol is a Dyson symbol and $F_1(m;n)=M(-m,n).$
The following theorem shows the function $\f_k(m_1,  \ldots, m_k;n)$ has the mirror symmetry with respect to each $m_j$.

 \begin{thm}\label{mirror-sym}  For $n\geq 2$,  $k\geq 1$ and   $1\leq j\leq k$,
  we have
\begin{align}
F_k(m_1,\ldots,m_j,\ldots,m_k;n)=F_k(m_1,\ldots,-m_j,\ldots,m_k;n).
\end{align}
\end{thm}

\pf   The above identity is trivial for $m_j=0$. We now assume that $m_j>0$. Let $H_k(m_1, \ldots,m_k;n)$ denote the set of $k$-marked Dyson symbols of $n$
counted by $F_k(m_1, \ldots,$\break $m_k;n)$. We aim to build a bijection $\Lambda$ between the set $H_k(m_1,\ldots,m_j,\ldots,m_k;n)$ and the set $H_k(m_1,\ldots,-m_j,\ldots,m_k;n)$.

Let
   \[\eta=\left(\begin{array}{ccccccccccccc}
\alpha^{(k)},&&\alpha^{(k-1)},&&\ldots,&&\alpha^{(j)},&&\ldots, &&&\alpha^{(1)}\\[2pt]
&p_{k-1},&&p_{k-2},&\cdots&p_{j}&&&&\cdots &p_1\\[2pt]
\beta^{(k)},&&\beta^{(k-1)},&&\ldots,&&\beta^{(j)},&&\ldots,&&&\beta^{(1)}
\end{array}
\right)\]
be a $k$-marked Dyson symbol  in $H_k(m_1,\ldots,m_j,\ldots,m_k;n)$.
To define the map $\Lambda$, we need to construct a new $j$-th vector $(\bar{\alpha}^{(j)},\bar{\beta}^{(j)})$   from  $( \alpha ^{(j)}, \beta ^{(j)})$. There are four cases.

 \noindent Case 1:   $1\leq j\leq k-1$. Set $\bar{\alpha}^{(j)}=\beta^{(j)}$ and $\bar{\beta}^{(j)}=\alpha^{(j)}$.

\noindent Case 2:  $j=k$ and $\ell(\alpha^{(k)})=1$. In this case, we have $\alpha_1^{(k)}=p_{k-1}$ and $\beta^{(k)}=\emptyset$.
 Set $\bar{\alpha}^{(k)}=\emptyset$ and $\bar{\beta}^{(k)}=\alpha^{(k)}.$

\noindent Case 3:  $j=k$, $\ell(\alpha^{(k)})\geq 2$ and $\ell(\beta^{(k)})\neq 1$.  Let $t=\beta_1^{(k)}-\beta_2^{(k)}$. Set
\[\bar{\alpha}^{(k)}=(\beta_1^{(k)}-t,\ \beta_2^{(k)},\ \ldots,\ \beta_\ell^{(k)}) \quad \text{and} \quad \bar{\beta}^{(k)}=(\alpha_1^{(k)}+t,\ \alpha_2^{(k)},\ \ldots,\ \alpha_\ell^{(k)}).\]

\noindent Case 4: $j=k$, $\ell(\alpha^{(k)})\geq 2$ and $\ell(\beta^{(k)})=1$.   Let $t=\beta_1^{(k)}-p_{k-1}$.  Set
\[\bar{\alpha}^{(k)}=(\beta_1^{(k)}-t) \quad \text{and} \quad \bar{\beta}^{(k)}=(\alpha_1^{(k)}+t,\ \alpha_2^{(k)},\ \ldots,\ \alpha_\ell^{(k)}).\]
From the above construction, it can be checked that
 \[\ell(\bar{\alpha}^{(j)})-\ell(\bar{\beta}^{(j)})
 =-(\ell(\alpha^{(j)})-\ell(\beta^{(j)})).\]

Then $\Lambda(\eta)$ is defined as
\[ \left(\begin{array}{ccccccccccccc}
\alpha^{(k)},&&\alpha^{(k-1)},&&\ldots,&&\bar{\alpha}^{(j)},&&\ldots, &&&\alpha^{(1)}\\[2pt]
&p_{k-1},&&p_{k-2},&\cdots&p_{j}&&&&\cdots &p_1\\[2pt]
\beta^{(k)},&&\beta^{(k-1)},&&\ldots,&&\bar{\beta}^{(j)},&&\ldots,&&&\beta^{(1)}
\end{array}
\right).\]
Hence $\Lambda(\eta)$ is a $k$-marked Dyson symbol in $H_k(m_1,\ldots,-m_j,\ldots,m_k;n)$. Furthermore, it can be seen that the above process is reversible. Thus $\Lambda$ is a bijection. \qed

We are now ready to prove
 Theorem \ref{main}, which says that  the number of $k$-marked Dyson symbols of $n$
 can be expressed in terms of the number of Dyson symbols of $n$.
This theorem is needed in  the combinatorial interpretation of $\mu_{2k}(n)$ given in
Theorem  \ref{main-w-4}. By  Theorem \ref{mirror-sym}, we see that
  Theorem \ref{main} can be deduced from the following formula.

\begin{thm}\label{main-p-f}
For $n\geq 2$ and $m_1,m_2,\ldots,m_k\geq 0$, we have
\begin{equation}\label{re-kmd-d}
\f_k(m_1,\ldots,m_k;n)=\sum_{t_1,\ldots,\,t_{k-1}=0}^{+\infty}
F_1\left(\sum_{i=1}^k m_i+2\sum_{i=1}^{k-1}t_i+k-1;n\right).
\end{equation}
\end{thm}

To prove the above theorem, we introduce the structure of strict  $k$-marked  Dyson symbols. Recall that a  strict  bipartition of $n$  is   a pair of partitions $(\alpha,\beta)$  such that  $\alpha_{i}>\beta_{i}$
 for $i=1,2,\ldots,\ell(\beta)$ and $|\alpha|+|\beta|=n.$
 Note that for a strick bipartition $(\alpha, \beta)$ we have
 $\ell(\alpha)\geq \ell (\beta)$.
 For example,
 \[\left(\begin{array}{cccccc}
 3&3&2&2&1\\
 2&1&1&1
\end{array}
\right)\]
 is a   strict   bipartition.

 Strict  bipartitions are the building blocks of
 strict  $k$-marked Dyson symbols. For $k\geq 2$, let
  \[\eta=\left(\begin{array}{cccccccc}
\alpha^{(k)},&&\alpha^{(k-1)},&&\ldots,&&\alpha^{(1)}\\[2pt]
&p_{k-1},&&p_{k-2},&\cdots&p_1\\[2pt]
\beta^{(k)},&&\beta^{(k-1)},&&\ldots,&&\beta^{(1)}
\end{array}
\right)\]
 be a $k$-marked Dyson symbols of $n$. If $(\alpha^{(i)},\beta^{(i)})$ is a strict bipartition for any  $1\leq i< k$,   we  say that  $\eta$   a  strict $k$-marked Dyson symbol of $n$.

 Notice that there is no balanced part in a strict  bipartition. Consequently, if $\eta$ is a strict $k$-marked Dyson symbol, then the $i$-th balanced number $b_i(\eta)$ of $\eta$ equals zero  for $1\leq i<k$.  To prove Theorem \ref{main-p-f}, we define a function $\f^{s}_k(m_1,\ldots,m_k;n)$ as the number of strict  $k$-marked Dyson symbols of $n$ with the $i$-th crank equal to $m_i$ for $1\leq i \leq k$ and define a function $\f_k(m_1,\ldots,m_k,t_1,\ldots,t_{k-1};n)$ as the number of    $k$-marked Dyson symbols of $n$ with the $i$-th crank equal to $m_i$ for $1\leq i \leq k$ and the $i$-th balance number equal to $t_i$ for $1\leq i \leq k-1$.
The relation stated in Theorem  \ref{main-p-f} can be established  via two steps as stated in the following two theorems.

\begin{thm}\label{lem-1}
For  $n\geq 2$,  $k\geq 2$,  $m_1,m_2,\ldots,m_k\geq 0$ and $t_1,t_2,\ldots,t_{k-1}\geq 0$, we have
\begin{equation}
\f_k(m_1,\ldots,m_k,t_1,\ldots,t_{k-1};n)=
 \f_k^{s}(m_1+2t_1,\ldots,m_{k-1}+2t_{k-1},m_k;n).
\end{equation}
\end{thm}

\begin{thm}\label{lem-2}
For   $n\geq 2$,  $k\geq 2$ and $m_1,m_2,\ldots,m_k\geq 0$, we have
\begin{equation}
\f_k^{s}(m_1,\ldots,m_k;n)=\f_1\left(\sum_{i=1}^km_i+k-1;n\right).
\end{equation}
\end{thm}

To prove Theorem \ref{lem-1}, we need  a  bijection    in  \cite[Theorem 2.4]{Ji-2011}. Let $P(r;n)$ denote the set of   pairs of partitions $(\alpha,\beta)$ of $n$  where there are  $r$ balanced parts
and  $\ell(\alpha)-\ell(\beta)\geq 0$, and let $Q(r;n)$ denote the set of strict bipartitions $(\bar{\alpha},\bar{\beta})$ of $n$  with $\ell(\bar{\alpha})-\ell(\bar{\beta})\geq r.$ Given two positive integers   $n$ and    $r$,  there is a   bijection $\psi$ between $P(r;n)$ and $Q(2r;n)$. Furthermore, the bijection $\psi$ possesses
 the following properties. For $(\alpha,\beta) \in P(r;n)$,
  let $(\bar{\alpha},\bar{\beta})=\psi(\alpha,\beta)$. Then we have
\begin{align}
&\bar{\alpha}_1=\max\{\alpha_1,\beta_1\},\quad   \bar{\alpha}_\ell= \alpha_\ell, \quad \text{and} \quad \bar{\beta}_\ell \geq \beta_\ell. \\[5pt]
&\ell(\bar{\alpha})=\ell(\alpha)+r  \quad \text{and} \quad \ell(\bar{\beta})=\ell(\beta)-r.
\end{align}

We next give a proof of  Theorem \ref{lem-1} by using the bijection $\psi$.

\noindent{\it Proof of Theorem \ref{lem-1}.} Let $P_k(m_1,\ldots,m_k,t_1,t_2,\ldots,t_{k-1};n)$ denote the set of $k$-marked Dyson symbols of $n$ with the $i$-th crank equal to
$m_i$ and the $i$-th balanced number
 equal to $t_i$,  and let $Q_k(m_1,\ldots,m_k;n)$ denote the set of strict $k$-marked  Dyson symbols of $n$ with  the $i$-th crank equal to $m_i$. We proceed to define a bijection $\Omega$ between  $P_k(m_1,\ldots,m_k,t_1,t_2,\ldots,t_{k-1};n)$ and  $Q_k(m_1+2t_1,\ldots,m_{k-1}+2t_{k-1}, m_k;n)$.

Let
\[\eta=\left(\begin{array}{cccccccc}
\alpha^{(k)},&&\alpha^{(k-1)},&&\ldots,&&\alpha^{(1)}\\[2pt]
&p_{k-1},&&p_{k-2},&\cdots&p_1\\[2pt]
\beta^{(k)},&&\beta^{(k-1)},&&\ldots,&&\beta^{(1)}
\end{array}
\right)\]
 be a $k$-marked Dyson symbol in $P_k(m_1,\ldots,m_k,t_1,t_2,\ldots,t_{k-1};n)$. For $1\leq i<k$,  we   apply the bijection   $\psi$ described above  to $(\alpha^{(i)},\beta^{(i)})$   to get  a pair of partitions $(\bar{\alpha}^{(i)},\bar{\beta}^{(i)})$. From
the properties of the bijection $\psi$, we see that$(\bar{\alpha}^{(i)},\bar{\beta}^{(i)})$ is a strict bipartition and
\begin{equation}\label{alph1}
\bar{\alpha}^{(i)}_1=\max\{\alpha^{(i)}_1,\beta^{(i)}_1\},\quad   \bar{\alpha}^{(i)}_\ell= \alpha^{(i)}_\ell,  \quad \bar{\beta}^{(i)}_\ell \geq \beta^{(i)}_\ell
\end{equation}
and
\begin{equation}\label{llalpha}
\ell(\bar{\alpha}^{(i)})=\ell(\alpha^{(i)})+t_i, \quad  \ell(\bar{\beta}^{(i)})=\ell(\beta^{(i)})-t_i.
\end{equation}
Then $\Omega(\eta)$ is defined to be
\[\left(\begin{array}{cccccccc}
\alpha^{(k)},&&\bar{\alpha}^{(k-1)},&&\ldots,&&\bar{\alpha}^{(1)}\\[2pt]
&p_{k-1},&&p_{k-2},&\cdots&p_1\\[2pt]
\beta^{(k)},&&\bar{\beta}^{(k-1)},&&\ldots,&&\bar{\beta}^{(1)}
\end{array}
\right). \]
By \eqref{alph1}, we see that  that for $1\leq i< k-1$,  each part of $\bar{\alpha}^{(i)}$ and $\bar{\beta}^{(i)}$ is between $p_{i-1}$ and $p_i$, namely,
    \[p_i\geq \bar{\alpha}^{(i)}_{1}\geq \bar{\alpha}^{(i)}_{2}\geq \cdots \geq \bar{\alpha}^{(i)}_{\ell}\geq p_{i-1} \quad \text{and} \quad  p_i\geq \bar{\beta}^{(i)}_{1}\geq \bar{\beta}^{(i)}_{2}\geq \cdots \geq \bar{\beta}^{(i)}_{\ell}\geq p_{i-1}.\]
It is also clear from \eqref{llalpha} that the $i$-th crank of $\Omega(\eta)$ is equal to $m_i+2t_i$ for $1\leq i<k$ and the $k$-th crank of $\Omega(\eta)$ is equal to $m_k$. Using   (\ref{llalpha}) again, we get
\[l(\Omega{(\eta)})=\sum_{i=1}^{k-1}\ell(\bar{\alpha}^{(i)})+\ell(\alpha^k)
=\sum_{i=1}^k(\ell(\alpha^{(i)})+t_i)=\sum_{i=1}^k\ell(\alpha^{(i)})+D
=l(\eta)+D
\]
and
\[s(\Omega{(\eta)})=\sum_{i=1}^{k-1}\ell(\bar{\beta}^{(i)})+\ell(\beta^k)
=\sum_{i=1}^k(\ell(\beta^{(i)})-t_i)=\sum_{i=1}^k\ell(\alpha^{(i)})-D=s(\eta)-D.
\]
Thus the weight of $\Omega{(\eta)}$ is equal to
\begin{align*}
 &\quad \sum_{i=1}^k(|\bar{\alpha}^{(i)}|+|\bar{\beta}^{(i)}|)
+\sum_{i=1}^{k-1}p_i+(l(\Omega{(\eta)})+k-1)\cdot s(\Omega{(\eta)})\\[3pt]
&\quad \quad =\sum_{i=1}^k(|\alpha^{(i)}|+|\beta^{(i)}|)+\sum_{i=1}^{k-1}p_i
+(l(\eta)+k-1+D)\cdot (s(\eta)-D),
\end{align*}
which is in accordance with the definition of $|\eta|$. So   $\Omega{(\eta)}$ is   in $Q_k(m_1+2t_1,\ldots,m_{k-1}+2t_{k-1}, m_k;n)$.
Since $\psi$ is a bijection, it is readily verified that
   $\Omega$ is  also a bijection, and hence the proof is complete. \qed

We  now turn to the proof of Theorem \ref{lem-2}.

\noindent{\it Proof of Theorem \ref{lem-2}.} Recall that $Q_k(m_1,\ldots,m_k;n)$ denotes the set of strict $k$-marked  Dyson symbols of $n$ with  the $i$-th crank equal to $m_i$ and  $H_1(m;n)$ denotes the set of Dyson symbols of $n$ with crank $m$.
To establish a bijection $\Phi$ between   $Q_k(m_1,\ldots,m_k;n)$ and  $H_1(m_1+\cdots+m_k+k-1;n)$, let
\[\eta=\left(\begin{array}{cccccccc}
\alpha^{(k)},&& \alpha^{(k-1)},&&\ldots,&&\alpha^{(1)}\\[2pt]
&p_{k-1},&&p_{k-2},&\cdots&p_1\\[2pt]
\beta^{(k)},&&\beta^{(k-1)},&&\ldots,&&\beta^{(1)}
\end{array}
\right)\]
be a strict $k$-marked  Dyson symbol in $Q_k(m_1,\ldots,m_k;n)$. Let $\alpha$ be the partition consisting of all parts of
$\alpha^{(1)} ,\alpha^{(2)},\ldots,\alpha^{(k)}$ together with   $p_1,\ldots,p_{k-1}$,
 and let $\beta$ be the partition
consisting of all parts of $\beta^{(1)},\beta^{(2)},\ldots,\beta^{(k)}$.  Then  $\Phi(\eta)$ is defined to be $(\alpha,\beta).$  From the definition of $k$-marked Dyson symbols, we see that $(\alpha,\beta)$ is
a Dyson symbol. It is also  easily seen that
\begin{equation}\label{dif-lengt}
 \ell(\alpha)=l(\eta)+k-1, \quad \ell(\beta)=s(\eta)
\end{equation}
and
\begin{equation}\label{dif-weight}
|\alpha|=\sum_{i=1}^k|\alpha^{(i)}|+\sum_{i=1}^{k-1}p_i, \quad  |\beta|=\sum_{i=1}^k|\beta^{(i)}|.
\end{equation}
It follows from \eqref{dif-lengt} that
\[\ell(\alpha)-\ell(\beta)=\sum_{i=1}^k m_i+k-1.\]
Combining  \eqref{dif-lengt} and \eqref{dif-weight}, we deduce that the weight of $(\alpha,\beta)$  equals
\begin{align*}
&|\alpha|+|\beta|+\ell(\alpha)  \ell(\beta)=\sum_{i=1}^k|\alpha^{(i)}|+\sum_{i=1}^{k-1}p_i
+\sum_{i=1}^k|\beta^{(i)}|+(l(\eta)+k-1) s(\eta)  =|\eta|.
\end{align*}
This proves that $(\alpha,\beta)$ is a Dyson symbol in $H_1(m_1+\cdots+m_k+k-1;n)$.

We next describe the reverse map  of $\Phi$.  Let
\[\left(\begin{array}{l}
\alpha\\
\beta
\end{array}\right)=\left(\begin{array}{cccc}
\alpha_1&\alpha_2&\ldots&\alpha_\ell\\[2pt]
\beta_1&\beta_2&\ldots&\beta_\ell
\end{array}
\right)\] be a Dyson symbol in $H_1(m_1+\cdots+m_k+k-1;n)$.
We proceed to show that a strict $k$-marked Dyson symbol $\eta$ can be recovered from
the Dyson symbol $( \alpha, \beta)$.

First,  we see that the $k$-th vector $(\alpha^{(k)},\beta^{(k)})$ of $\eta$ and $p_{k-1}$ can be recovered from $(\alpha,\beta)$.   Let $j_k$ be largest nonnegative integer such that
$\beta_{j_k}\geq \alpha_{m_k+j_k+1}$,  that is, for any $i\geq j_k+1$,
we have $\beta_{i}< {\alpha}_{m_k+i+1}$. Define
\[\left(\begin{array}{l}
\alpha^{(k)}\\
\beta^{(k)}
\end{array}\right)=\left(\begin{array}{cccc}
 {\alpha}_1& {\alpha}_2&\ldots& {\alpha}_{m_k+j_k}\\[2pt]
\beta_1&\beta_2&\ldots&\beta_{j_k}
\end{array}
\right)\quad \text{and} \quad p_{k-1}=\alpha_{m_k+j_k+1}.\]   Obviously,
$\ell(\alpha^{(k)})-\ell(\beta^{(k)})=m_k.$

To recover   $(\alpha^{(k-1)},\beta^{(k-1)})$ and $p_{k-1}$, we   let
\[\left(\begin{array}{l}
 {\alpha}'\\
\beta'
\end{array}\right)=\left(\begin{array}{cccc}
{\alpha}_{m_k+j_k+2}& {\alpha}_{m_k+j_k+3}&\ldots& {\alpha}_{\ell}\\[2pt]
\beta_{j_k+1}&\beta_{j_k+2}&\ldots&\beta_{\ell}
\end{array}
\right).\]
By the choice of $j_k$, we find that ${\alpha}_{m_k+j_k+i+1}>\beta_{j_k+i}$ for
any $i$, in other words,  $\alpha'_{i}>\beta'_i$. Consequently,  $( {\alpha}',\beta')$ is a strict bipartition.  Then $(\alpha^{(k-1)},\beta^{(k-1)})$  and $p_{k-1}$ can be constructed from $(\alpha',\beta')$. Let $j_{k-1}$
be the largest nonnegative integer such that $\beta'_{j_{k-1}}\geq
 {\alpha}'_{m_{k-1}+j_{k-1}+1}$. Define
\[\left(\begin{array}{l}
\alpha^{(k-1)}\\
\beta^{(k-1)}
\end{array}\right)=\left(\begin{array}{cccc}
 {\alpha}'_1& {\alpha}'_2&\ldots& {\alpha}'_{m_{k-1}+j_{k-1}}\\[2pt]
\beta'_1&\beta'_2&\ldots&\beta'_{j_{k-1}}
\end{array}
\right)\quad \text{and}  \quad p_{k-2}={\alpha}'_{m_{k-1}+j_{k-1}+1}. \]
Now we have $\ell(\alpha^{(k-1)})-\ell(\beta^{(k-1)})=m_{k-1}$.   Since $({\alpha}',\beta')$ is a strict bipartition, we deduce  that
$(\alpha^{(k-1)},\beta^{(k-1)})$ is a  strict bipartition.

The above procedure can be repeatedly used to determine
$(\alpha^{(k-2)},\beta^{(k-2)}), p_{k-3}, \ldots,$\break $(\alpha^{(2)},\beta^{(2)}),p_1,(\alpha^{(1)},\beta^{(1)})$.
The $k$-marked Dyson symbol $\eta$ can be defined as
\[ \left(\begin{array}{cccccccc}
\alpha^{(k)},&& \alpha^{(k-1)},&&\ldots,&&\alpha^{(1)}\\[2pt]
&p_{k-1},&&p_{k-2},&\cdots&p_1\\[2pt]
\beta^{(k)},&&\beta^{(k-1)},&&\ldots,&&\beta^{(1)}
\end{array}
\right).\]
It can be checked that
$\eta$ is  a strict $k$-marked  Dyson symbol in $Q_k(m_1,\ldots,m_k;n)$.
Moreover, it can be seen that $\Phi(\eta)=(\alpha,\beta)$, that is, $\Phi$
is indeed a bijection.  This completes the proof. \qed

Here is an example to illustrate the reverse map $\Phi^{-1}$.
 Assume that $m_1=1,m_2=1,m_3=0$, and
\[\left(\begin{array}{l}
 {\alpha}\\
\beta
\end{array}\right) =\left(\begin{array}{ccccccccccc}
6&6&3&3&3&3&2&2&1&1&1\\
5&5&4&2&1&1&1
\end{array}
\right), \]
which a Dyson symbol of $127$, that is, $(\alpha,\beta)\in H_1(4;127)$.
From  $({\alpha},
\beta)$,   we get
 \[
 \left(\begin{array}{l}
\alpha^{(3)}\\
\beta^{(3)}
\end{array}\right)=\left(\begin{array}{cccc}
6&6&3\\
5&5&4
\end{array}
\right),\quad p_2=3, \quad \left(\begin{array}{l}
 {\alpha}'\\
\beta'
\end{array}\right)=\left(\begin{array}{cccccccc}
 3&3&2&2&1&1&1\\
 2&1&1&1&&&
\end{array}
\right).\]
Based on  $( {\alpha}', \beta')$, we get
\[
 \left(\begin{array}{l}
\alpha^{(2)}\\
\beta^{(2)}
\end{array}\right)=\left(\begin{array}{cccccc}
 3&3&2&2&1\\
2&1&1&1
\end{array}
\right), \quad p_2=1, \quad \left(\begin{array}{l}
\alpha^{(1)}\\
\beta^{(1)}
\end{array}\right)=\left(\begin{array}{cc}
1\\&
\end{array}
\right).
\]
Finally, we obtain
\begin{equation*}\label{example}
\eta=\left(\begin{array}{cccccccccccccccc}
(6\ 6\ 3)&& (3\  3\  2\  2\  1)&&( 1)\\
&3&&1\\
(5\  5\  4 )&&(2\ 1\ 1\ 1)&&
\end{array}
\right).\end{equation*}
It can be checked that $\eta\in Q_3(1,1,0;127)$.

\section{A combinatorial interpretation of $\mu_{2k}(n)$}

In this section, we use Theorem \ref{main} to give a combinatorial interpretation of $\mu_{2k}(n)$  in terms of $k$-marked Dyson symbols.

\begin{thm}\label{main-w-4}For $k\geq 1$ and $n\geq 2$,   $\mu_{2k}(n)$ is equal to the number of $(k+1)$-marked Dyson symbols of $n$.
\end{thm}

\pf  By definition of $F_k(m_1,\ldots,m_k;n)$, the assertion of the theorem
can be stated as follows
\begin{align}\label{w-3-eq-1}
&\sum_{m_1, \ldots,m_{k+1}=-\infty}^{\infty}
\f_{k+1}(m_1,\ldots,m_{k+1};n)=\mu_{2k}(n).
\end{align}
Using Theorem \ref{main}, we see that
the left-hand side of \eqref{w-3-eq-1} equals
\begin{align}\label{3.2}
&\sum_{m_1,m_2,\ldots,m_{k+1}=-\infty}^{\infty}
\f_{k+1}(m_1,\ldots,m_{k+1};n)\nonumber \\[3pt]
&\hskip 2cm =\sum_{m_1,m_2,\ldots,m_{k+1}=-\infty}^{\infty}
\sum_{t_1,\ldots,t_{k}=0}^{\infty}
F_1\left(\sum_{i=1}^{k+1}|m_i|+2\sum_{i=1}^{k}t_i+k;n\right).
\end{align}
Given $k$ and $n$, let $c_k(j)$ denote the number of integer solutions to the equation
\[
|m_1|+\cdots+|m_{k+1}|+2t_1+\cdots+2t_k=j
\]
in $m_1, m_2, \ldots, m_{k+1}$ and $t_1, t_2, \ldots, t_k$
subject to the further condition that $t_1, t_2, \ldots, t_k$ are
nonnegative.
It can be shown that generating function of $c_k(j)$ is equal to
\begin{align*}
\sum_{j=0}^{\infty}c_k(j)q^j
%&=\left(\frac{1+q}{1-q}\right)^{k+1}\left(\frac{1}{1-q^2}\right)^k
%\nonumber\\[3pt]
&=\frac{1+q}{(1-q)^{2k+1}}, \label{gf-ck}
\end{align*}
so that
\begin{align*}
c_k(j)={2k+j \choose 2k}+{2k+j-1 \choose 2k}.
\end{align*}
Substituting $j$ by $m-k$, we get
$$c_k(m-k)={m+k-1 \choose 2k}+{m+k \choose 2k}.$$
Thus \eqref{3.2} simplifies to
\begin{align*}
&\sum_{m_1,m_2,\ldots,m_{k+1}=-\infty}^{\infty}
\f_{k+1}(m_1,\ldots,m_{k+1};n)\nonumber \\[3pt]
&\hskip 2cm =\sum_{m=1}^{\infty}\left[{m+k-1 \choose 2k}+{m+k \choose 2k}\right]F_1(m;n).
\end{align*}
Using Corollary \ref{core-dyson} and noting that $M(-m,n)=M(m,n)$, we conclude that
\begin{align*}
&\sum_{m_1,m_2,\ldots,m_{k+1}=-\infty}^{\infty}
\f_{k+1}(m_1,\ldots,m_{k+1};n)\nonumber \\[3pt]
&\hskip 2cm =\sum_{m=1}^{\infty}\left[{m+k-1 \choose 2k}+{m+k \choose 2k}\right]M(m,n),
\end{align*}
which equals $\mu_{2k}(n)$, as claimed. \qed

For example, for $n = 5$ and $k = 1$, we have $\mu_2(5)=35$, and
there are 35 $2$-marked Dyson symbols of $5$ as listed in the following table.
\[\begin{array}{cccc}
&
\left(\begin{array}{lllll}\\
&&&&1\\
(1&1&1&1)& \end{array}\right)
&\left(\begin{array}{ll} (1)\\
 &1\\
 (1)\end{array}\right)
&\left(\begin{array}{lll}
 (2&2)&\\
 &&1\\
&&\end{array}\right) \\[30pt]
&\left(\begin{array}{ll}&(1)\\
2&\\
&( {2}) \end{array}\right)
&\left(\begin{array}{llll}
&&&(1)\\
&&1\\
 (1 &1)&&( {1})\end{array}\right)
&\left(\begin{array}{lll}
&(1&1)\\
1\\
&( {1}& {1})\end{array}\right)\\[30pt]
&\left(\begin{array}{lllll}
(1&1&1&1)&\\
&&&&1\\
&&&& \end{array}\right)
&\left(\begin{array}{lll}\\
&&1\\
 (2 &2) \end{array}\right)
&\left(\begin{array}{llll}
(1&1) &&(1)\\
 && 1\\
&&&( {1})\end{array}\right)\\[30pt]
&
\left(\begin{array}{lllll}
&(1&1&1&1)\\
1\\
&&&& \end{array}\right)
&\left(\begin{array}{lllll}
(1&1)&&(1&1)\\
&&1\\
&&&& \end{array}\right)
&\left(\begin{array}{lllll}
&&&&(1)\\
&&&1\\
(1&1&1)&& \end{array}\right)\\[30pt]
\end{array}
\]
\[
\begin{array}{cccc}
&\left(\begin{array}{lll}
(2)&&(1)\\
 &2\\
 && \end{array}\right)
&\left(\begin{array}{lllll}
&&(1& 1&1)\\
&1 \\
(1)&\end{array}\right)
&\left(\begin{array}{lllll}
(1&1&1)&&(1)\\
&&&1&\\
&&&& \end{array}\right)\\[30pt]
&\left(\begin{array}{llll}
 (1)&&(1&1)\\
 &1 &\\
&&( {1})&\end{array}\right)
&\left(\begin{array}{lllll}
(1)&&(1&1&1)\\
&1\\
&&&& \end{array}\right)
&\left(\begin{array}{lllll}
&&&(1&1)\\
&&1\\
(1&1)&& \end{array}\right)\\[30pt]
&\left(\begin{array}{lll}
&(2&1)\\
 2 \\
&&\end{array}\right)
&\left(\begin{array}{llll}
&(1&1&1)\\
1\\
&( {1})\end{array}\right)
&\left(\begin{array}{llll}
&&(1&1)\\
&1\\
(1)&&( {1})\end{array}\right)\\[30pt]
&\left(\begin{array}{lll}
&&(1)\\
&2\\
(2)  \end{array}\right)
&\left(\begin{array}{llll}
(1)&&( {1})\\
&1\\
&&(1&1)\end{array}\right)
&\left(\begin{array}{llll}
&( {1}) \\
1&\\
&(1&1&1)\end{array}\right)\\[30pt]
&\left(\begin{array}{lllll}\\
1&&&\\
&(1&1&1&1)\end{array}\right)
&\left(\begin{array}{lllll}(1)\\
&1&\\
&&(1&1&1) \end{array}\right)
&\left(\begin{array}{lllll}\\
&&1&\\
(1&1)&&(1&1) \end{array}\right)\\[30pt]
&\left(\begin{array}{lll}\\
 2 &\\
&(2&1)\end{array}\right)
&\left(\begin{array}{lllll}
(1&1)\\
&&1&\\
&&&(1&1) \end{array}\right)
&\left(\begin{array}{lllll}\\
&&&1&\\
(1&1&1)&&(1) \end{array}\right)\\[30pt]
&\left(\begin{array}{lll}(2)\\
 &2\\
 &&(1) \end{array}\right)
&\left(\begin{array}{lllll}\\
&1&\\
(1)&&(1&1&1) \end{array}\right)
&\left(\begin{array}{lllll}
(1&1&1)\\
&&&1\\
&&&&(1) \end{array}\right)\\[30pt]
&\left(\begin{array}{llll}
&&( {1})\\
&1\\
(1)&&(1&1)\end{array}\right)
&\left(\begin{array}{lll}\\
&2\\
(2)&&(1) \end{array}\right)
\end{array}
\]

\section{ Congruences for $\mu_{2k}(n)$ }

In this section, we introduce  the full crank of a $k$-marked Dyson symbol. We   show that there exist an infinite  family of congruences for the full crank function of $k$-marked Dyson symbols.

To define the full crank of a $k$-marked Dyson symbol $\eta$, denoted $FC(\eta)$, we
recall that $c_k(\eta)$ denotes the $k$-th  crank of $\eta$,  $l(\eta)$ denotes the large length of $\eta$ and $s(\eta)$ denotes the small  length of $\eta$ and  $D$ denotes the balanced number of $\eta$.  Then  $FC(\eta)$ is given by
 \[FC(\eta)=\begin{cases}
 l(\eta)-s(\eta)+2D+k-1, & \text{if} \quad  c_k(\eta)>0,  \\[3pt]
-(l(\eta)-s(\eta)+2D+k-1), & \text{if}  \quad c_k(\eta)\leq 0.
 \end{cases}\]
 It is clear that for $k=1$, the full crank of a $1$-marked Dyson symbol reduces to the crank of a  Dyson symbol.

  Analogous to the full rank function for a $k$-marked Durfee symbol defined by Andrews \cite{Andrews-07-a}, we define the full crank function  $\c_k(i,t;n)$  as the number of $k$-marked Dyson symbols of $n$ with the full crank congruent to $i$ modulo $t$.
 The following theorem gives an infinite family of congruences
 of the full crank function.

\begin{thm}\label{conguren-main} For fixed prime $p\geq 5$   and positive integers $r$ and $k\leq (p+1)/2$.    Then there exist infinitely many non-nested arithmetic progressions $An+B$ such that for each $0\leq i\leq p^r-1$,
\[\c_k(i, p^r; An+B) \equiv 0 \pmod{p^r}.\]
\end{thm}

Since
\[\mu_{2k}(n)=\sum_{i=0}^{p^r-1}\c_{k+1}(i,p^r;n),\]
 Theorem \ref{conguren-main} implies the following congruences for $\mu_{2k}(n)$.

\begin{thm}For  fixed prime $p\geq 5$,  positive integers $r$ and $k\leq (p-1)/2$. Then there exists infinitely many non-nested arithmetic progressions $An+B$ such that
\[\mu_{2k}(An+B)\equiv 0 \pmod{p^r}.\]
\end{thm}

To prove Theorem \ref{conguren-main}, let $\c_k(m;n)$  denote   the number of $k$-marked Dyson symbols of $n$ with the full crank equal to $m$. In this notation,  we have the following relation.

\begin{thm}\label{thm4.3}For $n\geq 2$, $k\geq 1$ and integer $m$,
\begin{equation}\label{NCk&Ck}
 \c_k(m;n)={m+k-2 \choose 2k-2}M(m,n).
\end{equation}
\end{thm}

\pf Recall that $\f_k(m_1,\ldots,m_k,t_1,\ldots,t_{k-1};n)$ is the number of    $k$-marked Dyson symbols of $n$ such that for $1\leq i \leq k$, the $i$-th crank equal to $m_i$  and the $i$-th balance number equal to $t_i$. By the definition of $\c_k(m,n)$, we see that if $m\geq 1$, then we have
\begin{equation}\label{4.2}
\c_k(m;n)=\sum F_k(m_1,m_2,\ldots,m_k,t_1,t_2,\ldots,t_{k-1};n),
\end{equation}
where the summation ranges over  all integer solutions to the equation
\begin{equation}\label{4.3}
|m_1|+\cdots+|m_{k}|+2t_1+\cdots+2t_{k-1}=m-k+1
\end{equation}
in $m_1, m_2, \ldots, m_{k}$ and $t_1, t_2, \ldots, t_{k-1}$
subject to the further condition that $m_k$ is positive  and $t_1, t_2, \ldots, t_{k-1}$ are nonnegative.

Combining Theorem  \ref{lem-1} and Theorem  \ref{lem-2}, we find that
\begin{equation}\label{4.4}
F_k(m_1,m_2,\ldots,m_k,t_1,t_2,\ldots,t_{k-1};n)=F_1\left(\sum_{i=1}^{k}|m_i|+2\sum_{i=1}^{k-1}t_i+k-1;n\right).
\end{equation}
Substituting  \eqref{4.4} into \eqref{4.2}, we get
\begin{equation}\label{labelck}
\c_k(m;n)=\sum F_1\left(\sum_{i=1}^{k}|m_i|+2\sum_{i=1}^{k-1}t_i+k-1;n\right),
\end{equation}
where the summation ranges over all  solutions to the equation \eqref{4.3}. Let $\bar{c}_k(m-k+1)$ denote the number of   integer solutions to the equation \eqref{4.3}. It is not difficult to verify that
\[
\bar{c}_k(m-k+1)= {m+k-2 \choose 2k-2}.
\]
Thus, \eqref{labelck} simplifies to
\begin{align*}
 \c_k(m;n)={m+k-2 \choose 2k-2}F_1(m;n).
\end{align*}
Using Corollary \ref{core-dyson} and noting that $M(-m,n)=M(m,n)$, we conclude that
 \begin{align*}
 \c_k(m;n)={m+k-2 \choose 2k-2}M(m,n),
\end{align*}
as required. Similarly,
it can be shown that relation \eqref{NCk&Ck} also holds for $m\leq 0$.  \qed

 Let $M(i,t;n)$ denote the number of partitions of $n$ with the crank congruent to $i$ modulo $t$. The following congruences for $M(i,t;n)$ given by  Mahlburg \cite{Mahlburg-2005} will be used in the proof of Theorem \ref{conguren-main}.

\begin{thm}[Mahlburg]\label{Mahl-thm}
For fixed prime $p\geq 5$  and positive integers $\tau$ and $r$,  there are infinitely many non-nested
arithmetic progressions $An+B$ such that for each $0\leq m\leq p^r-1$,
\[
M(m,p^{r};An+B)\equiv 0\pmod{p^{\tau}}.
\]
\end{thm}

We are now ready to complete the proof of  Theorem \ref{conguren-main} by
using Theorems \ref{thm4.3} and  \ref{Mahl-thm}.

\noindent{\it Proof of Theorem \ref{conguren-main}.} For     $0\leq i \leq p^r-1$,
by the definition of $\c_k(i,p^r;n)$, we have
\begin{equation}\label{4.6}
 \c_k(i,p^r;n)=\sum_{t=-\infty}^{+\infty}\c_k(p^rt+i;n).
\end{equation}
Replacing $m$ by $p^rt+i$ in \eqref{NCk&Ck}, we get
\begin{equation}\label{4.7}
\c_k(p^rt+i;n)={p^rt+i+k-2 \choose 2k-2}M(p^rt+i,n).
\end{equation}
Substituting \eqref{4.7} into \eqref{4.6}, we find that
\begin{equation}\label{NCk&Ck-2}
 \c_k(i,p^r;n)=\sum_{t=-\infty}^{+\infty}{p^rt+i+k-2 \choose 2k-2}M(p^rt+i,n).
\end{equation}
Since $p$ is a prime and  $k\leq (p+1)/2$, we see that
 $(2k-2)!$ is not divisible by $p$. It follows that
\[{p^rt+i+k-2 \choose 2k-2}\equiv {i+k-2 \choose 2k-2} \pmod{p^r}.\]
Thus \eqref{NCk&Ck-2} implies that
\begin{align*}\label{NCk&Ck-2}
 \c_k(i,p^r;n)&\equiv \sum_{t=-\infty}^{+\infty}{i+k-2 \choose 2k-2}M(p^rt+i,n) \pmod{p^{r}}\\[3pt]
 &={i+k-2 \choose 2k-2}M(i,p^r;n).
\end{align*}
Setting $\tau=r$ in Theorem \ref{Mahl-thm}, we deduce that   there are infinitely many non-nested
arithmetic progressions $An+B$ such that for every $0\leq i\leq p^r-1$
\[
M(i,p^{r};An+B)\equiv 0\pmod{p^{r}}.
\]
Consequently,  there are infinitely many non-nested
arithmetic progressions $An+B$ such that for every $0\leq m\leq p^r-1$
\[
NC_k(i,p^{r};An+B)\equiv 0\pmod{p^{r}},
\]
and hence the proof is complete. \qed

\noindent{\bf Acknowledgments.} This work was supported by the 973 Project and the National Science Foundation of China.

\end{document}